\documentclass[11pt]{article} 

\usepackage[utf8]{inputenc} 

\usepackage{geometry} 
\usepackage{graphicx} 

\usepackage{amssymb}

\usepackage{booktabs} 
\usepackage{array} 
\usepackage{paralist} 
\usepackage{verbatim} 
\usepackage{subfig} 

\usepackage{fancyhdr} 
\usepackage{sectsty}
\allsectionsfont{\sffamily\mdseries\upshape} 

\usepackage[nottoc,notlof,notlot]{tocbibind} 
\usepackage[titles,subfigure]{tocloft} 

\newtheorem{Thm}{Theorem}
\newtheorem{Prop}{Proposition}
\newtheorem{Lem}{Lemma}

\newtheorem{Ex}{Example}

\newcommand{\bpr}{\noindent \textbf{Proof}: ~}

\newcommand{\eprsk}{~$\blacksquare$\medskip}

\title{On $1$-cocycles induced by a positive definite function on a locally compact abelian group}
\author{Jordan Franks\footnote{Partially supported by a ThinkSwiss Research Scholarship.}\hspace{3mm}and Alain Valette}

\begin{document}
\maketitle
\begin{abstract} For $\varphi$ a normalized positive definite function on a locally compact abelian group $G$, we consider on the one hand the unitary representation $\pi_\varphi$ associated to $\varphi$ by the GNS construction, on the other hand the probability measure $\mu_\varphi$ on the Pontryagin dual $\hat{G}$ provided by Bochner's theorem. We give necessary and sufficient conditions for the vanishing of 1-cohomology $H^1(G,\pi_\varphi)$ and reduced 1-cohomology $\overline{H}^1(G,\pi_\varphi)$. For example, $\overline{H}^1(G,\pi_\varphi)=0$ if and only if either $Hom(G,\mathbb{C})=0$ or $\mu_\varphi(1_G)=0$, where $1_G$ is the trivial character of $G$.
\end{abstract}

\section{Introduction}

The Gel'fand-Naimark-Segal construction (see \cite{BHV08}) provides a correspondence  between positive definite functions $\varphi$ on a locally compact group $G$ and cyclic representations $\pi_\varphi$ on Hilbert space.  This allows  one to construct a dictionary between the functional-analytic and algebro-geometric pictures of $\varphi$ and $\pi_\varphi$. For example, $\varphi$ is an extreme point in the cone $\mathcal{P}(G)$ of positive definite functions on $G$ if and only if $\pi_\varphi$ is an irreducible representation; or, there exists a constant $a>0$ such that $\varphi-a$ is again positive definite if and only if $\pi_\varphi$ has nonzero fixed vectors (see \cite{Dix69}).  

In view of their importance for rigidity questions and Kazhdan's property (T)\footnote{Recall Shalom's results, see Theorems 0.2 and 6.1 in \cite{Sha00}: for a compactly generated group $G$, the group $G$ has property (T), if and only if $\overline{H}^1(G,\pi)=0$ for every unitary representation $\pi$ of $G$, if and only if $H^1(G,\sigma)=0$ for every unitary irreducible representation of $G$.}, it is natural to try to fit 1-cohomology and reduced 1-cohomology of $\pi_\varphi$ in that dictionary. This is the question we address in this paper, assuming $G$ to be a locally compact abelian group. Indeed, in this case we enjoy Bochner's Theorem (see \cite{Fol95}): $\varphi$  is the Fourier transform of a positive Borel measure $\mu_\varphi$ on the Pontryagin dual $\hat{G}$.  Without relying on the cohomological machinery available in the literature (see \cite{Gui80,BHV08}), we achieve by completely elementary means the results of this paper, namely, that the existence of nontrivial 1-cohomology is determined by two factors: the existence of non-trivial homomorphisms to $\mathbb{C}$, and more importantly, the behavior of $\mu_\varphi$ near the trivial character $1_G.$  

\section{Statement of results}

For $G$ a locally compact group and $\pi$ a unitary representation of $G$ on a Hilbert space ${\cal H}_\pi$, recall that the space of {\it 1-cocycles} for $\pi$ is
$$Z^1(G,\pi)=\{b:G\rightarrow{\cal H}_\pi: b\;\mbox{continuous},\,b(gh)=\pi(g)b(h)+b(g)\,\mbox{for all}\,g,h\in G\}.$$
The space of {\it 1-coboundaries} for $\pi$ is:
$$B^1(G,\pi)=\{b\in Z^1(G,\pi):\exists v\in{\cal H}_\pi\,\mbox{such that}\,b(g)=\pi(g)v-v\,\mbox{for every}\,g\in G\}.$$
The {\it 1-cohomology} of $\pi$ is then the quotient
$$H^1(G,\pi)=Z^1(G,\pi)/B^1(G,\pi).$$
Endow $Z^1(G,\pi)$ with the topology of uniform convergence on compact subsets of $G$. The {\it reduced 1-cohomology} of $\pi$ is the quotient of the space of 1-cocycles by the closure of the space of 1-coboundaries, i.e.
$$\overline{H}^1(G,\pi)=Z^1(G,\pi)/\overline{B^1(G,\pi)}.$$

From now on, let $G$ be a locally compact {\it abelian} group, $\varphi$ a positive definite function on $G$.  Excluding the zero function, we may without loss of generality take $\varphi$ to be normalized ($\varphi(e)=1$). Let $\mu_\varphi$ be the probability measure on the Pontryagin dual $\hat{G}$ provided by Bochner's theorem, i.e. $\varphi(x)=\int_{\widehat{G} }\xi(x)\,d\mu_\varphi(\xi)$ for $x\in G$.  Let $(\pi_\varphi,{\cal H}_\varphi,\xi_\varphi)$ be the cyclic representation of $G$ associated to $\varphi$ through the GNS construction, so that the cyclic vector $\xi_\varphi\in{\cal H}_\varphi$ satisfies $\langle \pi_\varphi(x)\xi_\varphi|\xi_\varphi\rangle=\varphi(x)$. Let also $\rho_\varphi$ be the representation of $G$ on $L^2(\hat{G},\mu_\varphi)$ given by $(\rho_\varphi(x)f)(\xi)=\xi(x)f(\xi)\;(\xi\in\hat{G},f\in L^2(\hat{G},\mu_\varphi)$.

If $\lambda_G$ denotes the regular representation of $G$ and $h\in L^2(G)$, then $\widehat{(\lambda_G(x)h)}(\xi)=\overline{\xi(x)}\hat{h}(\xi)$.  From Plancherel's Theorem, it follows that the composition of the Fourier transform with conjugation is a unitary equivalence between the regular representation on $L^2(G)$ and the unitary representation defined by $\Big(f(\xi)\mapsto \xi(x) f(\xi)\Big)$ on $L^2(\hat{G})$.  This, together with Bochner's Theorem, intuits our introduction of $\rho_\varphi$, as well as the following proposition which we prove in section 3.     

\begin{Prop}\label{equiv} The representations $\pi_\varphi$ and $\rho_\varphi$ are unitarily equivalent.
\end{Prop}

We assume from now on that $\varphi$ is not the constant function $1$, so that $\mu_\varphi$ is not the Dirac mass at the trivial character $1_G$ of $G$. This is still equivalent to $\mu_\varphi(1_G)<1$. Let $\mu_\varphi^\perp$ be the probability measure on $\hat{G}$ defined by $\mu_\varphi=\mu_\varphi(1_G)\delta_{1_G}+(1-\mu_\varphi(1_G))\mu_\varphi^\perp$. 

Let $\pi_\varphi^0$ be the (trivial) subrepresentation of $\pi_\varphi$ on the subspace ${\cal H}^0_\varphi$ of $\pi_\varphi$-fixed vectors, and $\pi_\varphi^\perp$ be the subrepresentation on the orthogonal complement, so that $\pi_\varphi =\pi_\varphi^0\oplus\pi_\varphi^\perp$. A simple computation in the $\rho_\varphi$-picture shows that ${\cal H}^0_\varphi\neq 0$ if and only if $\mu_\varphi(1_G)>0$, and in this case ${\cal H}^0_\varphi =\mathbb{C}\delta_{1_G}$. Moreover the map $L^2(\hat{G}\backslash\{1_G\},\mu_\varphi^\perp)\rightarrow L^2(\hat{G},\mu_\varphi):f\mapsto\frac{f}{\sqrt{1-\mu_\varphi(1_G)}}$ is isometric and identifies $\pi_\varphi^\perp$ with the restriction of $\rho_\varphi$ to $L^2(\hat{G}\backslash\{1_G\},\mu_\varphi^\perp)$.

Our main result is:

\begin{Thm}\label{coh} Let $\varphi$ be a nonconstant, normalized positive definite function on a locally compact abelian group $G$. 

\begin{enumerate}
\item[1)] Consider the following statements:
\begin{enumerate}
  \item[i)] $H^1(G,\pi_\varphi)= 0$;
  \item[ii)] Both of the following properties are satisfied:
  \begin{enumerate}
  \item[a)] $\mu_\varphi(1_G)= 0$ or $Hom(G,\mathbb{C})= 0$;
  \item[b)] $1_G\notin supp (\mu_\varphi^\perp)$.
\end{enumerate}
Then $(ii)\Rightarrow(i)$, and the converse holds if $G$ is $\sigma$-compact.
\end{enumerate}

\item[2)] The following are equivalent:
\begin{enumerate}
  \item[i)] $\overline{H}^1(G,\pi_\varphi)= 0$;
  \item[ii)] $\mu_\varphi(1_G)= 0$ or $Hom(G,\mathbb{C})= 0$.
  \end{enumerate}
\end{enumerate}
\end{Thm}

This result will be proved in section 4. It is essentially equivalent to Theorem 4 in \cite{Gui72}, but we emphasize the fact that our proof is direct and based on explicit construction of cocycles and coboundaries.

\section{Proof of Proposition \ref{equiv}}

\begin{Lem}\label{cyclic} The constant function $1\in L^2(\hat{G},\mu_\varphi)$ is a cyclic vector for $\rho_\varphi$.
\end{Lem}

\bpr For $f\in L^1(G)$, consider the operator $\rho_\varphi(f)=\int_G f(x)\rho_\varphi(x)\,dx$; then $(\rho_\varphi(f).1)(\xi)=\int_G f(x)\xi(x)\,dx=\hat{f}(\xi)$. 

Denote by $C_0(\hat{G})$ the space of continuous functions vanishing at infinity on $\hat{G}$, and recall that $\hat{f}\in C_0(\hat{G})$ (the Riemann-Lebesgue Lemma). It is classical that the map $L^1(G)\rightarrow C_0(\hat{G})$ is a continuous algebra homomorphism with dense image (a consequence of Stone-Weierstrass). Now compose this homomorphism with the continuous inclusion $C_0(\hat{G})\rightarrow L^2(\hat{G},\mu_\varphi): h\mapsto h.1$.  Since continuous functions with compact support are dense in $L^2(\hat{G},\mu_\varphi)$, this inclusion has dense image. Since the map $L^1(G)\rightarrow L^2(\hat{G},\mu_\varphi):f\mapsto\rho_\varphi(f).1$ is the composite of the previous maps, it has dense image, meaning that $1$ is cyclic for $\rho_\varphi$.
\eprsk

Observe that $\langle \rho_\varphi(x).1|1\rangle = \int_{\widehat{G}}\xi(x)\,d\mu_\varphi(\xi)=\varphi(x)$, so Proposition \ref{equiv} follows from Lemma \ref{cyclic} and the uniqueness statement of the GNS construction.
\eprsk

\section{Proof of Theorem \ref{coh}}

Since $\pi_\varphi=\pi_\varphi^0\oplus\pi_\varphi^\perp$, we have $H^1(G,\pi_\varphi)=H^1(G,\pi_\varphi^0)\oplus H^1(G,\pi_\varphi^\perp)$ and analogously for $\overline{H}^1$. As $B^1(G,\pi_\varphi^0)=0$, we see that $H^1(G,\pi_\varphi^0)= 0$ if and only if $\overline{H}^1(G,\pi_\varphi^0)= 0$, if and only if either $\mu_\varphi(1_G)= 0$ or $Hom(G, \mathbb{C})= 0$: this proves the implications $(i)\Rightarrow (ii)(a)$ in part 1 of Theorem \ref{coh}, and $(i)\Rightarrow (ii)$ in part 2 of Theorem \ref{coh}; moreover, it reduces the main result to:

\begin{Thm}\label{perp} Let $\varphi$ be a nonconstant, normalized positive definite function on a locally compact abelian group $G$.
\begin{enumerate}
\item[1)] If $1_G \notin supp(\mu_\varphi^\perp)$, then $H^1(G,\mu_\varphi^\perp)=0$. The converse holds if $G$ is $\sigma$-compact.
\item[2)] $\overline{H}^1(G,\pi_\varphi^\perp)=0$.
\end{enumerate}
\end{Thm}

\begin{Ex} In Part 1 of Theorems \ref{coh} and \ref{perp}, the converse implications are {\it false} when $G$ is not assumed to be $\sigma$-compact. Indeed, let $G$ be an uncountable abelian group with the discrete topology, and take $\varphi=\delta_1$. Then $\pi_\varphi$ is the left regular representation $\lambda_G$ on $\ell^2(G)$, while $\mu_\varphi=\mu_\varphi^\perp$ is the Haar measure on the compact group $\hat{G}$. Since $\mu_\varphi$ has full support, in particular $1_G$ lies in its support. On the other hand $H^1(G,\lambda_G)=0$ by Proposition 4.13 in \cite{CTV08}
\end{Ex} 

To prove the implication ``$\Rightarrow$'' in part 1 of Theorem \ref{perp}, we will need:

\begin{Lem}\label{proba} Let $F$ be a closed subset of $\hat{G}$, with $1_G\notin F$.
\begin{enumerate}
\item[a)] There exists a regular Borel probability measure $\nu_0$ on $G$ such that the Fourier transform $\widehat{\nu_0}$ vanishes on $F$.
\item[b)] For every $\varepsilon>0$, there exists a compactly supported regular Borel probability measure $\nu$ on $G$ such that $|\hat{\nu}|<\varepsilon$ on $F$.
\end{enumerate}
\end{Lem}

\bpr $(a)$ See section $1.5.2$ in \cite{Rud62}.

$(b)$ Let $\nu_0$ be a probability measure on $G$ as in $(a)$. Let $\delta$ be a number $0<\delta<1$, to be determined later. Let $C$ be a compact subset of $G$ such that $\nu_0(C)>1-\delta$. Let $\nu$ be the probability measure on $G$ defined by $\nu(B)=\frac{\nu_0(B\cap C)}{\nu_0(C)}$, for every Borel subset $B\subseteq G$. By taking $\delta$ small enough, the total variation distance $|\nu_0-\nu|(G)$ between $\nu_0$ and $\nu$ can be made arbitrarily small. For any finite signed measure $\mu$ on $G$, we have the classical inequality $|\int_G f(x)\,d\mu(x)|\leq\|f\|_\infty |\mu|(G)$; applied to $\mu=\nu_0-\nu$ and $f(x)=\xi(x)$ with $\xi\in\hat{G}$, it gives $|\widehat{\nu_0}(\xi)-\hat{\nu}(\xi)|\leq |\nu_0-\nu|(G)$, so that $\|\widehat{\nu_0}-\widehat{\nu}\|_\infty<\varepsilon$ for $\delta$ small enough.
\eprsk

\medskip
\noindent
{\bf Proof of ``$\Rightarrow$'' in part 1 of Theorem \ref{perp}:} We assume that $1_G$ is not in the support of $\mu_\varphi^\perp$, and prove that $H^1(G,\pi_\varphi^\perp)=0$. Let $b\in Z^1(G,\pi_\varphi^\perp)$ be a 1-cocycle. Expanding $b(xy)=b(yx)$ using the cocycle relation, we get:
$$(1-\pi_\varphi^\perp(x))b(y)=(1-\pi_\varphi^\perp(y))b(x)\;\;\;(x,y\in G).$$
In the realization of $\pi_\varphi^\perp$ on $L^2(\hat{G},\mu_\varphi^\perp)$, this gives:
\begin{equation}\label{coc}
(1-\xi(x))b(y)(\xi)=(1-\xi(y))b(x)(\xi)
\end{equation}
almost everywhere in $\xi$ (w.r.t. $\mu_\varphi^\perp$). By Lemma \ref{proba}, we can find a compactly supported probability measure $\nu$ on $G$ such that $|1-\hat{\nu}|\geq\frac{1}{2}$ on $supp(\mu_\varphi^\perp)$. Define an element $v\in L^2(\hat{G},\mu_\varphi^\perp)$ by $v:=\int_G b(y)\,d\nu(y)$: since $b$ is continuous and $\nu$ is compactly supported, the integral exists (in the weak sense) in $L^2(\hat{G},\mu_\varphi^\perp)$. Integrating (\ref{coc}) w.r.t. $\nu$ in the variable $y$, we get:
\begin{equation}\label{cocintegr}
(1-\xi(x))v(\xi)=(1-\hat{\nu}(\xi))b(x)(\xi)
\end{equation}
almost everywhere in $\xi$. Since $|1-\hat{\nu}|\geq\frac{1}{2}$ on $supp(\mu_\varphi^\perp)$, the function $w(\xi):=\frac{v(\xi)}{1-\hat{\nu}(\xi)}$ belongs to $L^2(\hat{G},\mu_\varphi^\perp)$, and by (\ref{cocintegr}) its coboundary is exactly $b$.
\eprsk

\begin{Lem}\label{nonzerocoh} Let $H$ be a locally compact group. Let $(\sigma_n)_{n\geq 1}$ be a sequence of unitary representations of $H$ without nonzero fixed vectors, with $\sigma_n$ acting on a Hilbert space ${\cal H}_n$. Assume that, for each $n\geq 1$, there exists a unit vector $\eta_n\in{\cal H}_n$ such that the series $\sum_{n=1}^\infty\|\sigma_n(x)\eta_n-\eta_n\|^2$ converges uniformly on compact subsets of $H$. Set $\sigma=\oplus_{n=1}^\infty \sigma_n$. Then $b(x):=\oplus_{n\geq 1}^\infty (\sigma_n(x)\eta_n-\eta_n)$ defines a nonzero element in $H^1(G,\sigma)$.
\end{Lem}

\bpr By assumption  $b(x)$ belongs to $\oplus_{n=1}^\infty{\cal H}_n$ and the map $H\rightarrow \oplus_{n=1}^\infty{\cal H}_n:x\mapsto b(x)$ is continuous. Let $\eta\in\prod_{n=1}^\infty {\cal H}_n$ be defined as $\eta=(\eta_n)_{n\geq 1}$. Since $b$ is the formal coboundary of $\eta$, we have $b\in Z^1(H,\sigma)$. To prove that $b$ is not a coboundary, it suffices to show that the associated affine action $\alpha(x)v=\sigma(x)v+b(x)$ on $\oplus_{n=1}^\infty{\cal H}_n$ has no fixed point. But $\alpha(x)v=v$ translates into $\sigma_n(x)(v_n+\eta_n)=v_n+\eta_n$ for every $x\in H$ and $n\geq 1$. Since $\sigma_n$ has no nonzero fixed vector, we have $v_n+\eta_n=0$ so $\|v_n\|=1$, which contradicts $\sum_{n=1}^\infty \|v_n\|^2 <+\infty$.
\eprsk

\medskip
\noindent
{\bf Proof of ``$\Leftarrow$'' in part 1 of Theorem \ref{perp}, assuming $G$ to be $\sigma$-compact:} Let $(K_n)_{n\geq 0}$ be an increasing sequence of compact subsets of $G$, with $G=\cup_{n=1}^\infty K_n$, and $K_0=\{1\}$. Define a basis $(U_k)_{k\geq 0}$ of open neighborhoods of $1_G$ in $\hat{G}$ by $U_k=\{\xi\in\hat{G}: \max_{g\in K_k} |\xi(g)-1|<2^{-k}\}$ (observe that $U_0=\hat{G}$). Define a sequence $(k_n)_{n\geq 0}$ inductively by $k_0=0$ and $k_n=\min\{k: k>k_{n-1},\,\mu_\varphi^\perp(U_k)<\mu_\varphi^\perp(U_{k_{n-1}})\}$ for $n\ge 1$ (since $\mu_\varphi^\perp\{1_G\}=0$ and $1_G$ is in the support of $\mu_\varphi^\perp$, this is well-defined). Set then $C_n:=U_{k_{n}}\backslash U_{k_{n+1}}$ for $n\ge 1$, and let ${\cal H}_n$ be the space of functions in $L^2(\hat{G},\mu_\varphi^\perp)$ which are $\mu_\varphi^\perp$-almost everywhere 0 on $\hat{G}\backslash C_n$. Then ${\cal H}_n$ is a closed, $\rho_\varphi$-invariant subspace of $L^2(\hat{G},\mu_\varphi^\perp)$.  Denote by $\sigma_n$ the restriction of $\rho_\varphi$ to ${\cal H}_n$, so that $L^2(\hat{G},\mu_\varphi^\perp)=\oplus_{n=0}^\infty{\cal H}_n$ and $\rho_\varphi=\oplus_{n=0}^\infty \sigma_n$. Let $\eta_n=\frac{{\bf 1}_{C_n}}{\sqrt{\mu_\varphi^\perp(C_n)}}$ be the normalized characteristic function of $C_n$. To appeal to Lemma \ref{nonzerocoh}, we still have to check that $x\mapsto \sum_{n=1}^\infty\|\sigma_n(x)\eta_n-\eta_n\|^2$ converges uniformly on every compact subset $K$ of $G$. Clearly we may assume $K=K_\ell$. For $n\geq\ell$ and $x\in K_\ell$ and $\xi\in C_n$, we have $|\xi(x)-1|< 2^{-k_n}$, hence:
$$\max_{x\in K_\ell} \sum_{n=\ell}^\infty\|\sigma(x)\eta_n-\eta_n\|^2=\max_{x\in K_\ell}\sum_{n=\ell}^\infty \frac{1}{\mu_\varphi^\perp(C_n)}\int_{C_n} |\xi(x)-1|^2\,d\mu_\varphi^\perp(\xi)\leq \sum_{n=\ell}^\infty 4^{-k_n}\leq \sum_{n=0}^\infty 4^{-n}=\frac{4}{3}$$
and
$$\max_{x\in K_\ell} \sum_{n=0}^\infty\|\sigma(x)\eta_n-\eta_n\|^2\leq(\max_{x\in K_\ell} \sum_{n=0}^{\ell-1}\|\sigma(x)\eta_n-\eta_n\|^2) + \frac{4}{3}\leq 4\ell + \frac{4}{3}<+\infty.$$
So the result follows from Lemma \ref{nonzerocoh}.
\eprsk

\medskip
\noindent
{\bf Proof of part 2 of Theorem \ref{perp}:} Let $b\in Z^1(G,\pi_\varphi^\perp)$ be a 1-cocycle.  We must show that $b$ is a limit of 1-coboundaries (uniformly on compact subsets of $G$). Since $\mu_\varphi^\perp(1_G)=0$, by regularity of $\mu_\varphi^\perp$, we may find a decreasing sequence of relatively compact open neighborhoods $(V_n)_{n\ge 1}$ of $1_G$, such that $\mu_\varphi^\perp(V_n)\rightarrow 0$ for $n\rightarrow\infty$. Set ${\cal H}_n:=\{f\in L^2(\hat{G},\mu_\varphi^\perp): f=0\;\mbox{a.e. on $V_n$}\}$; then ${\cal H}_n$ is a closed, $\rho_\varphi$-invariant subspace, and the sequence $({\cal H}_n)_{n\ge 1}$ is increasing with dense union in $L^2(\hat{G},\mu_\varphi^\perp)$. Let $\rho_n$ denote the restriction of $\rho_\varphi$ to ${\cal H}_n$, and $b_n$ be the projection of $b$ onto ${\cal H}_n$. Then $b_n\in Z^1(G,\rho_n)$, and $\lim_{n\rightarrow\infty} b_n=b$ (uniformly on compact subsets of $G$). But since $1_G$ does not belong to the closed subset $\hat{G}\backslash V_n$, by part 1 of Theorem \ref{perp} we have $H^1(G,\rho_n)=0$, so that $b_n$ is a coboundary. 
\eprsk


\newpage
Authors'addresses:

\medskip
J.F.: Mathematisches Institut, Universit\"at Bonn, Endenicher Allee 60, 53115 Bonn, Germany, E-mail: \textsf{jjfranks@outlook.com}.  

\medskip
A.V.: Institut de Math\'ematiques, Universit\'e de Neuch\^atel, Unimail, 11 Rue Emile Argand, CH-2000 Neuch\^atel, Switzerland. E-mail: \textsf{alain.valette@unine.ch}
\end{document}